\documentclass{amsart}

\usepackage[english,british]{babel}
\usepackage{amssymb,amsmath,mathtools,xcolor,graphicx,xspace,colortbl,ragged2e,rotating} %
\usepackage{amsfonts}  
\usepackage{amsmath}  
\usepackage{amsthm}  
\usepackage[activate={true,nocompatibility},final,tracking=true,kerning=true,spacing=true,factor=1100,stretch=10,shrink=10]{microtype}
\usepackage[pdftex,pdfborderstyle={/S/U/W 1},hyperfootnotes=false,hidelinks]{hyperref}
\usepackage{autonum}
\newcommand{\tmtextit}[1]{{\itshape{#1}}}
\newcommand{\tmstrong}[1]{\textbf{#1}}
\newtheorem {theorem}{Theorem}

 \newtheorem{proposition}[theorem]{Proposition}

\begin{document}
\title[A note on the paper ``The Cremona problem in
dimension 2'']{A note on the paper ``The Cremona problem in
dimension 2'' by Wolfgang Bartenwerfer
}
\author[S.~Brzostowski]{Szymon Brzostowski}
\address{Faculty~of~Mathematics and~Computer~Science\\
University of {\L}{\'o}d{\'z}\\
ul. Banacha 22, 90-238 {\L}{\'o}d{\'z}\\
Poland\\}
\email{szymon.brzostowski@wmii.uni.lodz.pl}

\author[T.~Krasi\'{n}ski]{Tadeusz Krasi\'{n}ski}
\address{Faculty~of~Mathematics and~Computer~Science\\
University of {\L}{\'o}d{\'z}\\
ul. Banacha 22, 90-238 {\L}{\'o}d{\'z}\\
Poland\\}
\email{tadeusz.krasinski@wmii.uni.lodz.pl}

\begin{abstract}


The paper titled ``Cremona problem in dimension 2'' by W.~Bar\-tenwerfer presented a flawed attempt at proving the Jacobian Conjecture. Our aim is to provide a thorough analysis of the author's approach, highlighting the errors that were made in the process.
\end{abstract}

\maketitle

\section{Introduction
}
In issue 1 of volume 119 of this journal the article ``The Cremona problem in
dimension 2'' by Wolfgang Bartenwerfer has appeared. In this article, the famous Jacobian Conjecture in the two-dimensional case was solved positively.


{\bigskip}
{\noindent}{\tmstrong{Jacobian Conjecture. }}\tmtextit{If a polynomial mapping $F :\mathbb{C}^{n} \rightarrow \mathbb{C}^{n}$ has a non-zero constant jacobian ($\det \ensuremath{\operatorname*{Jac}}F =const. \neq 0$), then $F$ is a polynomial automorphism. }{\bigskip}

\noindent This is one of the famous problems posed by Stephen Smale at the turn of the century \cite{S}. It has been studied by many mathematicians and several times (allegedly) solved. A lot of information on this subject, one may find in two monographs devoted to this topic, by the leading figure in this domain -- Arno van den Essen, with co-authors \cite{E1}, \cite{E2}. So, the article by Wolfgang Bartenwerfer in this journal was a surprise to the experts in this field. However, there are some weak points in this proof of the Jacobian Conjecture.

In the article we explain these obscure points, provide details of the reasonings, and point out the faults of the author. These observations should convince everyone that the proof couldn't be repaired following the idea presented in the mentioned paper. However, there are some advantages of the paper. Methods of extension of the field $\mathbb{C}$ to another algebraically closed, complete field endowed with a non-archimedean absolute value and transfering considerations to the Laurent series (at infinity) in two complex variables are valuable and may help to decide the Jacobian Conjecture. A consequence of our article is that the Jacobian Conjecture is still open (even for $n =2$). 

\section{
The idea of the ``proof'' along W. Bartenwerfer. 
}
The method of the proof is \textit{reductio ad absurdum}; the assumption that the Jacobian Conjecture is false leads to a contradiction. So, assume that $F : =(f ,g) :\mathbb{C}^{2} \rightarrow \mathbb{C}^{2}$ is a \textit{ jacobian pair,} i.e., $\det \ensuremath{\operatorname*{Jac}}F =1$ (such mappings are sometimes called \textit{Keller mappings}) and $F$ is not a polynomial automorphism. Since $(f ,g)$ is a jacobian pair, the Abhyankar theorem (Cor. 18.15 in \cite{A}) asserts that $f$ and $g$ have at most two points at infinity. So, by a linear change of coordinates in $\mathbb{C}^{2} ,$ we may assume about $g$
\begin{equation}g(X ,Y) =X^{n_{1}}Y^{n_{2}} +t .l .d. ,\quad n : =n_{1} +n_{2} >0 ,\quad n_{1} \geq 0 ,\; n_{2} \geq 0.
\end{equation}($t .l .d.$ means terms of lower degrees, hence $\deg g =n$). By another Abhyankar theorem   (Thm 19.4 in \cite{A}):

\begin{theorem}
The following conditions are equivalent:

\noindent 1.$\quad\forall _{(\widetilde{f} ,\widetilde{g}) \in \mathbb{C}[X ,Y]}(\det \ensuremath{\operatorname*{Jac}}(\widetilde{f} ,\widetilde{g}) =1 \Leftrightarrow (\widetilde{f} ,\widetilde{g}) \text{ is a polynomial automorphism})$

\noindent 2.$\quad\forall _{(\widetilde{f} ,\widetilde{g}) \in \mathbb{C}[X ,Y]}(\det \ensuremath{\operatorname*{Jac}}(\widetilde{f} ,\widetilde{g}) =1 \Leftrightarrow \widetilde{g} \text{ has only one point at infinity})$
\end{theorem}

\noindent we may assume that $g$ has exactly two points at infinity, i.e.,
\begin{equation}g(X ,Y) =X^{n_{1}}Y^{n_{2}} +t .l .d. ,\quad n_{1} >0 , \; n_{2} >0.
\end{equation}
Hence $n =n_{1} +n_{2} \geq 2.$ Simple algebraic considerations, using the \textit{jacobian condition} $\det \ensuremath{\operatorname*{Jac}}F =1 ,$  give that $f$ has  the leading term proportional to that of $g$
\begin{gather}f(X ,Y) =c_{1}X^{m_{1}}Y^{m_{2}} +t .l .d .,\quad m_{1} >0 , \; m_{2} >0 , \; c_{1} \neq 0 , \\
\frac{m_{1}}{n_{1}} =\frac{m_{2}}{n_{2}} . \nonumber \end{gather}
Hence
\begin{equation}
\frac{m}{n} =\frac{m_{1}}{n_{1}} =\frac{m_{2}}{n_{2}}, \label{5}
\end{equation}

\noindent where $m : =m_{1} +m_{2} =\deg f$. We may additionally assume that $m \neq n .$ In fact, if $m =n$ then by (\ref{5}) $m_{1} =n_{1} ,$ $m_{2} =n_{2}$ and then we may change the pair $(f ,g)$ to the pair $(f -\frac{1}{c_{1}}g ,g) ,$ which satisfies the required condition and still is a jacobian pair and is not a polynomial automorphism. From the above the author infers that \begin{equation}\left \vert m_{1} -m_{2}\right \vert  \neq \left \vert n_{1} -n_{2}\right \vert  .
\end{equation}This is the first fault of the author because this inference is not justified (it was communicated to us by Z. Jelonek and A. Parusi\'{n}ski). In fact, from (\ref{5}) we get
\begin{equation}\left \vert m_{1} -m_{2}\right \vert  =\frac{m}{n}\left \vert n_{1} -n_{2}\right \vert  .
\end{equation}
Since $\frac{m}{n} \neq 1$, this implies\begin{equation}\left \vert m_{1} -m_{2}\right \vert  \neq \left \vert n_{1} -n_{2}\right \vert \text{    or  
}\left \vert m_{1} -m_{2}\right \vert  =\vert n_{1} -n_{2}\vert  =0.\text{} \label{8}
\end{equation}The second case $m_{1} =m_{2} ,$ $n_{1} =n_{2}$ has not been considered further by the author. However, the next considerations will give $n_{1} \neq n_{2} ,$ which shows that the first case holds.  So, this is not a decisive error.

In order to use (\ref{5}) for simplifying the pair $(f ,g)$ (which is obvious to do when $\frac{m}{n} \in \mathbb{N})$, the author transfers the considerations to the ring ($\mathbb{C}$-algebra) of formal Laurent series at infinity in two variables $\mathbb{C}\mathbb{L}_{\infty }^{X ,Y}$ with coefficients in $\mathbb{C} .$ By definition
\begin{equation}\mathbb{C}\mathbb{L}_{\infty }^{X ,Y} : =\left\{\sum \limits _{i = -\infty }^{k}f_{i} :\quad \deg f_{i} =i , \; k \in \mathbb{Z}\right\} ,
\end{equation}
where  $f_{i}$  are Laurent polynomials in two variables with coefficients in $\mathbb{C}$, i.e., elements of $\mathbb{C}[X ,Y]_{XY}$ -- the localization of $\mathbb{C}[X ,Y]$ with respect to the multiplicative set $\{(XY)^{i} :i =0 ,1 ,\ldots \} .$ In fact $\mathbb{C}[X ,Y]_{XY} =\mathbb{C}[X ,X^{ -1} ,Y ,Y^{ -1}] .$ In the $\mathbb{C}$-algebra $\mathbb{C}\mathbb{L}_{\infty }^{X ,Y} ,$ for our specific element $g(X ,Y) ,$ we may speak of its inverse element and a $k$-th root of it, for appropriate $k .$ Namely, we have  
\begin{align}
g(X ,Y) &= X^{n_{1}}Y^{n_{2}}(1 +u(X ,Y)) ,\quad u \in \mathbb{C}\mathbb{L}_{\infty }^{X ,Y} , \; \deg u <0 , \nonumber  \\
g(X ,Y)^{ -1} &= X^{ -n_{1}}Y^{ -n_{2}}\frac{1}{1 +u(X ,Y)}\\
 &=X^{ -n_{1}}Y^{ -n_{2}}(1 -u(X ,Y) +u(X ,Y)^{2} - \cdots ) , \nonumber
 \end{align}
and for any $k ,$ $k\vert n_{1} ,k\vert n_{2}$ we may write a $k$-th root of $g(X ,Y)$ 
\begin{align}
g^{\frac{1}{k}}(X ,Y) =X^{\frac{n_{1}}{k}}Y^{\frac{n_{2}}{k}}\sqrt[k]{1+u(X,Y)} =X^{\frac{n_{1}}{k}}Y^{\frac{n_{2}}{k}}\sum \limits _{i =0}^{\infty }\left (\genfrac{}{}{0pt}{}{1/k}{i}\right )u(X ,Y)^{i}. \label{11}
\end{align}
In the latter case, we get in fact $k$ roots. They are
\begin{equation}g^{\frac{1}{k}}(X ,Y) ,\varepsilon g^{\frac{1}{k}}(X ,Y) ,\ldots  ,\varepsilon ^{k -1}g^{\frac{1}{k}}(X ,Y),
\end{equation} 

\noindent where $\varepsilon $ is a primitive $k$-th root of unity. Then the author, through a finite number of steps, using (\ref{5}), modifies (in $\mathbb{C}\mathbb{L}_{\infty }^{X ,Y}$) the pair $(f ,g)$ to another jacobian pair $(f -G ,g) \in (\mathbb{C}\mathbb{L}_{\infty }^{X ,Y})^{2} ,$ where $G$ ``depends'' on $g ,$ with the degree of $f -G$ as small as possible -- precisely degree $2 -n .$ However, it is not possible in some cases, which gives us (in these cases) the required contradiction. We will explain it in a while. The first step of the reduction of the degree of $f$ is simple: if we denote $\frac{l_{1}}{k_{1}} : =\frac{m_{1}}{n_{1}} =\frac{m_{2}}{n_{2}} ,$ $\ensuremath{\operatorname*{GCD}}(l_{1} ,k_{1}) =1 ,$ then
\begin{equation}
(g^{\frac{1}{k_{1}}})^{l_{1}}(X ,Y) =X^{m_{1}}Y^{m_{2}} +t .l .d.  \in \mathbb{C}\mathbb{L}_{\infty }^{X ,Y} .
\end{equation}
Hence
\begin{equation}\deg (f -c_{1}(g^{\frac{1}{k_{1}}})^{l_{1}}) <m_{1} +m_{2} =m
\end{equation}
and $(f -c_{1}(g^{\frac{1}{k_{1}}})^{l_{1}} ,g)$ is still a jacobian pair in $\mathbb{C}\mathbb{L}_{\infty }^{X ,Y} .$ Since $\ensuremath{\operatorname*{GCD}}(l_{1} ,k_{1}) =1 ,$ for all the $k_{1}$-th roots of unity $\varepsilon  \neq 1$ the expression $c_{1}(\varepsilon g^{\frac{1}{k_{1}}})^{l_{1}}$ does not cancel the leading term of $f .$ Hence  $\deg (f -c_{1}(\varepsilon g^{\frac{1}{k_{1}}})^{l_{1}}) =m_{1} +m_{2} .$ Notice that by the jacobian condition
\begin{equation}
	\deg (f -c_{1}(g^{\frac{1}{k_{1}}})^{l_{1}})  \geq 2 -n .
\end{equation}
If the inequality is sharp, again by the jacobian condition, we get that
$f -c_{1}(g^{\frac{1}{k_{1}}})^{l_{1}}$ has the leading term proportional to that of $g$
\begin{gather}
f -c_{1}(g^{\frac{1}{k_{1}}})^{l_{1}} =c_{2}X^{\widetilde{m}_{1}}Y^{\widetilde{m}_{2}} +t .l .d. , \text{ where }\widetilde{m}_{1} ,\widetilde{m}_{2} \in \mathbb{Z} , \; c_{2} \in \mathbb{C}^{ \ast },\\
\frac{\widetilde{m}_1}{n_1} = \frac{\widetilde{m}_2}{n_2}.
\end{gather}
 We continue this process until finally, after a finite number of steps, we get
  \begin{equation}
 \deg (f -c_{1}(g^{\frac{1}{k_{1}}})^{l_{1}} -c_{2}(g^{\frac{1}{k_{2}}})^{l_{2}} - . . . -c_{p}(g^{\frac{1}{k_{p}}})^{l_{p}}) =_{}2 -n, \label{17}
\end{equation}
\noindent where $k_{i} \in \mathbb{N}$, $l_{i} \in \mathbb{Z}$, $\ensuremath{\operatorname*{GCD}}(l_{i}, k_{i}) =1$, $k_{i}\vert n_{1}$, $k_{i}\vert n_{2}$ $(i =1 ,\ldots  ,p) $. Since $2 -n \leq 0 ,$ at least one of denominators $k_{i}$ is $ >1$ (otherwise we would have all $l_{i} >0 ,$ which implies $\deg (f -c_{1}(g^{\frac{1}{k_{1}}})^{l_{1}} -c_{2}(g^{\frac{1}{k_{2}}})^{l_{2}} - . . . -c_{p}(g^{\frac{1}{k_{p}}})^{l_{p}}) >0 ,$ which contradicts $2 -n \leq 0$). In particular, if $\ensuremath{\operatorname*{GCD}}(n_{1} ,n_{2}) =1$ we obtain a contradiction.
So, in the sequel we will assume $\ensuremath{\operatorname*{GCD}}(n_{1} ,n_{2}) >1$.
In case (\ref{17}) holds, again from the jacobian condition, the author infers (Proposition 1.3) that $n_{1} \neq n_{2}$ and the leading form of $f -c_{1}(g^{\frac{1}{k_{1}}})^{l_{1}} -c_{2}(g^{\frac{1}{k_{2}}})^{l_{2}} - . . . -c_{p}(g^{\frac{1}{k_{p}}})^{l_{p}}$ of degree $2 -n$ has at most two terms
\begin{equation}dX^{1 -n_{1}}Y^{1 -n_{2}} +c_{i ,j}X^{i}Y^{j} ,\enspace d \in \mathbb{C}^{ \ast } ,c_{i ,j} \in \mathbb{C} ,
\end{equation}
the second one provided there exist $i ,j \in \mathbb{Z} ,$ $i +j =2 -n ,$ $\frac{2 -n}{n} =\frac{i}{n_{1}} =\frac{j}{n_{2}} .$ In fact, the second term is impossible because from inequalities $0 <n_{i} <n$ and these conditions we easily get $n_{1} =n_{2} =n/2$ which contradicts $n_{1} \neq n_{2} .$ Hence
\begin{equation}
f -c_{1}(g^{\frac{1}{k_{1}}})^{l_{1}} -c_{2}(g^{\frac{1}{k_{2}}})^{l_{2}} - . . . -c_{p}(g^{\frac{1}{k_{p}}})^{l_{p}} =dX^{1 -n_{1}}Y^{1 -n_{2}} +t .l .d ., \enspace d \in \mathbb{C}^{ \ast } .
\end{equation}
If we denote $\widetilde{k} : =\ensuremath{\operatorname*{LCM}}(k_{1} ,\ldots  ,k_{p}) ,$ then $\widetilde{k} >1$ and we may represent the left hand side of the above formula as $f -W({g^{{1}/{{\widetilde{k}}}}})$,
 where $W(T)$ is a Laurent polynomial in one variable $T$. Precisely, $W(T) =c_{1}T^{l_{1}\widetilde{k}_{1}^{}} +\ldots  +c_{p}T^{l_{p}\widetilde{k}_{p}^{}} ,$ where $\widetilde{k}_{i} : =\widetilde{k}/k_{i} .$
 From the above procedure of reduction of the degree, it follows $W$ is unique and for each $\widetilde{k}$-th root of unity $\varepsilon  \neq 1$ we have $W(\varepsilon g^{{1}/{\widetilde{k}^{}}}) \neq W(g^{{1}/{\widetilde{k}}})$. In fact, if $\varepsilon ^{\widetilde{k}} =1 ,$ $\varepsilon  \neq 1 ,$ and  $W(\varepsilon g^{{1}/{\widetilde{k}}}) =W(g^{{1}/{\widetilde{k}}}) ,$ then we take the smallest $q ,$ $1 <q \leq \widetilde{k} ,$ such that $\varepsilon ^{q} =1.$ Then $s \cdot q =\widetilde{k}$ for some $s ,$ $1 \leq s <\widetilde{k} .$ We have $q\vert l_{1}\widetilde{k}_{1} ,\ldots  ,q\vert l_{p}\widetilde{k}_{p} .$ Hence $k_{1}\vert s ,\ldots  ,k_{p}\vert s .$ Since $\widetilde{k} =\ensuremath{\operatorname*{LCM}}(k_{1} ,\ldots  ,k_{p}) ,$ we get $\widetilde{k}\vert s$ which is impossible. Thus, again from the process of reduction, it follows that for such $\varepsilon _{} ,$ $\varepsilon  \neq 1$
 \begin{equation}
	\deg (f -W(\varepsilon g^{{1}/{\widetilde{k}}})) \neq 2 -n ,
\end{equation}and more precisely\begin{equation}f -W(\varepsilon g^{{1}/{\widetilde{k}}}) =eX^{n_{1}^{ \prime }}Y^{n_{2}^{ \prime }} +t .l .d .,{\enspace} n_{1}^{ \prime } +n_{2}^{ \prime } >2 -n, \; e \neq 0 .
\end{equation}
However, the author represents $W(g^{{1}/{\widetilde{k}}})$ as $\widetilde{W}(g^{{1}/{k}})$, where $k =\ensuremath{\operatorname*{GCD}}(n_{1} ,n_{2}) .$ Since $\widetilde{k}\vert k$ we have $\widetilde{W}(T) =W(T^{k/\widetilde{k}}) .$ Then we may only assert that there exists $\varepsilon _{0} ,$ $\varepsilon _{0}^{k} =1 ,$ such that $\widetilde{W}(\varepsilon _{0}g^{{1}/{k}}) \neq \widetilde{W}^{\vphantom{{A}^{A^a}}}(g^{{1}/{k}})$. 

After all these preparations, we must consider only pairs $(f ,g)$ satisfying the following \tmstrong{Conditions} (changing back the name $\widetilde{W}$ to $W)$:

\begin{enumerate}
\item\label{nu1} $n_{1} \neq n_{2} .$ This implies by (\ref{8}) that $\left \vert m_{1} -m_{2}\right \vert  \neq \left \vert n_{1} -n_{2}\right \vert \text{.}$

\item\label{nu2} $f -W(g^{\frac{1}{k}}) =dX^{1 -n_{1}}Y^{1 -n_{2}} +t .l .d . ,$  $k =\ensuremath{\operatorname*{GCD}}(n_{1} ,n_{2}) >1,$ $d \neq 0$.

\item\label{nu3} $f -W(\varepsilon _{0}g^{\frac{1}{k}}) =eX^{n_{1}^{ \prime }}Y^{n_{2}^{ \prime }} +t .l .d. ,$ $n_{1}^{ \prime } +n_{2}^{ \prime } >2 -n$,  $e \neq 0$, for some root of unity $\varepsilon _{0}$ of degree $k$.
\end{enumerate}

To get a contradiction under the above assumptions, the author extends the domain of solutions of the system $\left \{f(X ,Y) =\alpha  ,\; g(X ,Y) =\beta \right \}$ from complex numbers  to their extension containing the Laurent series at infinity in one variable. This is a natural idea as the behaviour of the mapping at infinity is crucial for the  Jacobian Conjecture. He extends $\mathbb{C}$ to a complete non-archimedean field $\mathbb{K}$ which may be defined as the last object in the chain of extensions
\begin{align}\mathbb{C} \subset \mathbb{C}[[t^{ -1}]] \subset \mathbb{C}((t^{ -1})) \subset \overline{\mathbb{C}((t^{ -1}))} \subset \widehat{\overline{\mathbb{C}((t^{ -1}))}} = :\mathbb{K} ,\end{align}where $\overline{K}$ means an algebraic closure of $K$ and $\widehat{R}$ -- a completion of $R$  with respect to the absolute value $\left \vert  \cdot \right \vert \text{ }$generated by the natural valuation $v(\varphi ) : = -\deg \varphi $ (then $\left \vert \varphi \right \vert  =e^{ -v(\varphi )})$ on $\mathbb{C}[[t^{ -1}]]$ and extended to the rest of the chain. As it is known, $\overline{\mathbb{C}((t^{ -1}))}$ is the ring of the formal Puiseux series at infinity. Notice that for $a \in \mathbb{C}^{ \ast }$ we have $v(a) =0$ and consequently $\left \vert a\right \vert  =1.$ We have an obvious extension $\mathbb{C}\mathbb{L}_{\infty }^{X ,Y} \subset \mathbb{K}\mathbb{L}_{\infty }^{X ,Y}$.

After this preparation, the author applies the non-archimedean absolute value in $\mathbb{K}$ to get all the above objects convergent. He considers the set \begin{equation}V =\{(x ,y) \in \mathbb{K}^{2} :1 <\varepsilon  \leq \left \vert x\right \vert  =\left \vert y\right \vert  \leq \varrho \}
\end{equation}for $\varepsilon  ,\varrho $ appropriately chosen, in which the above transformations lead to convergent series. So, we may assume
\begin{align}
f(X ,Y) &=X^{m_{1}}Y^{m_{2}}(1 +v(X ,Y)) ,\quad v \in \mathbb{K}\mathbb{L}_{\infty }^{X ,Y} , \; \deg v <0 , \\
g(X ,Y) &=X^{n_{1}}Y^{n_{2}}(1 +u(X ,Y)) ,\quad u \in \mathbb{K}\mathbb{L}_{\infty }^{X ,Y} , \; \deg u <0 , \\
g^{\frac{1}{k}}(X ,Y) &=X^{\frac{n_{1}}{k}}Y^{\frac{n_{2}}{k}}(1 +w(X ,Y)) ,\quad w \in \mathbb{K}\mathbb{L}_{\infty }^{X ,Y} , \; \deg w <0 , \label{g} \\
(f -W(g^{\frac{1}{k}}))&(X ,Y) \\&=dX^{1 -n_{1}}Y^{1 -n_{2}}(1 +z(X ,Y)) ,\quad z \in \mathbb{K}\mathbb{L}_{\infty }^{X ,Y} , \; \deg z <0 , \label{fw}
\end{align}
the Laurent series $v ,u ,w ,z$ are convergent in $V$ (in the metric induced by the non-archimedean absolute value in $\mathbb{K}$) and moreover $\left \vert u(x ,y)\right \vert  <1 ,\left \vert v(x ,y)\right \vert  <1 ,\left \vert w(x ,y)\right \vert  <1 ,\left \vert z(x ,y)\right \vert  <1$ for $(x ,y) \in V .$ It will be useful in the sequel that in fact coefficients of $v ,u ,w ,z$ belong to $\mathbb{C}$, i.e., $v ,u ,w ,z \in \mathbb{C}\mathbb{L}_{\infty }^{X ,Y} .$
The series  $f(X ,Y) ,$ $g(X ,Y) ,$ $g^{\frac{1}{k}}(X ,Y) ,$ $(f -W(g^{\frac{1}{k}}))(X ,Y)$ are treated as functions from $V$ into $\mathbb{K} .$ The author counts the cardinality of the fibres of the restriction\begin{equation}(f ,g)\vert _{V} :V \rightarrow (f ,g)(V)
\end{equation}
and claims that  the cardinality of each fibre is equal to $\left \vert n_{1} -n_{2}\right \vert  .$ The next author's steps, which lead to a contradiction, are as follows. If we interchange $f$ and $g$ in the above reasonings, then after shrinking $V$ we would get that the cardinality of each fibre of the same mapping $(f ,g)\vert _{V} :V \rightarrow (f ,g)(V)$ is equal to $\left \vert m_{1} -m_{2}\right \vert  .$ Hence $\left \vert m_{1} -m_{2}\right \vert  =\vert n_{1} -n_{2}\vert  ,$ which gives a contradiction with Condition \ref{nu1}. However, his computation of the cardinality of fibres of $(f ,g)\vert _{V} :V \rightarrow (f ,g)(V)$ contains a gap. His procedure to count this cardinality is as follows. Fix an arbitrary $(x_{0} ,y_{0}) \in V$ and put $\alpha _{0} : =f(x_{0} ,y_{0}) ,$ $\beta _{0} : =g(x_{0} ,y_{0})$ and $\widetilde{\beta }_{0} : =g^{\frac{1}{k}}(x_{0} ,y_{0}) .$ If $\left \vert x_{0}\right \vert  =\left \vert y_{0}\right \vert  = :r ,$ then from the basic properties of non-archimedean absolute value we get $\left \vert \alpha _{0}\right \vert  =r^{m} ,$ $\left \vert \beta _{0}\right \vert  =r^{n}$ , $\vert \widetilde{\beta} _{0}\vert  =r^{n/k} .$  So, we have to find solutions in $V$ of the system

\begin{equation}
  \left\{\!\!\!\begin{array}{l}
    f (X, Y) = \alpha_0\\
    g (X, Y) = \beta_0
  \end{array}\right.\!\!\!. \label{27}
\end{equation}
Notice that for every point $(x_{1} ,y_{1}) \in V$ with another absolute value $\left \vert x_{1}\right \vert  =\left \vert y_{1}\right \vert  =r_{1} \neq r$ we have $\left \vert f(x_{1} ,y_{1})\right \vert  =r_{1}^{m} \neq r^{m}$. Hence the fibre $(f ,g)\vert _{V}^{ -1}(\alpha _{0} ,\beta _{0})$ lies in the \textit{torus} $T_{r ,r} : =\{(x ,y) \in \mathbb{K}^{2} :\left \vert x\right \vert  =\left \vert y\right \vert  =r\} .$ So, it suffices to find solutions of system (\ref{27}) in $T_{r ,r} .$ But this set of solutions 
in $T_{r ,r}$  is obviously equal to the union of the sets of solutions in $T_{r ,r}$ of $k$ systems

\begin{equation}
  \left\{\!\!\!\begin{array}{l}
    f (X, Y) = \alpha_0\\
    \varepsilon^i g^{1 / k} (X, Y) = \widetilde{\beta}_0
  \end{array}\right. \label{19}
\end{equation}
$(i =0 , . . . ,k -1)$, where $\varepsilon $ is a $k$-th primitive root of $1.$ The author claims that each of these latter systems has the same number of solutions, equal to  $\vert n_{1} -n_{2}\vert /k ,$ and in consequence system (\ref{27}) has $\vert n_{1} -n_{2}\vert $ solutions. This is not justified because we have the following facts.

\begin{proposition} \label{prop2}
System (\ref{19}) for a fixed $i$ has $\vert n_{1} -n_{2}\vert /k$ solutions in $T_{r ,r}$ if and only if there exists at least one solution of this system in $T_{r ,r} .$
\end{proposition}

\begin{proof}
Assume system (\ref{19}) has a  solution $(x_{1} ,y_{1}) \in T_{r ,r} ,\left \vert x_{1}\right \vert  =\left \vert y_{1}\right \vert  =r .$ Notice that solutions of (\ref{19}) in $T_{r ,r}$ are the same as solutions of the slightly modified system
\begin{equation}
  \left\{\!\!\!\begin{array}{l}
    (f - W (g^{\frac{1}{k}}))  (X, Y) = \alpha_0 - W (\varepsilon^{- i} 
    \widetilde{\beta}_0)\\
    g^{1 / k} (X, Y) = \varepsilon^{- i}  \widetilde{\beta}_0
  \end{array}\right.\!\!\!. \label{ww}
\end{equation}
From forms (\ref{fw}), (\ref{g}) of $f -W(g^{\frac{1}{k}})$ and $g^{\frac{1}{k}}$ we get that the mapping $(f -W(g^{\frac{1}{k}}) ,g^{\frac{1}{k}})$ maps $T_{r ,r}$ into  $T_{r^{2 -n} ,r^{n/k}} .$ Moreover, by assumption on $(x_{1} ,y_{1})$, the point $(\alpha _{0} -W(\varepsilon ^{ -i}\widetilde{\beta }_{0}) ,\varepsilon ^{ -i}\widetilde{\beta }_{0})$ belongs to $T_{r^{2 -n} ,r^{n/k}}$ because

\begin{align}
\vert \alpha _{0} -W(\varepsilon ^{ -i}\widetilde{\beta }_{0}) \vert  &=\vert (f -W(g^{\frac{1}{k}}))(x_{1} ,y_{1})\vert  =\vert dx_{1}^{1 -n_{1}}y_{1}^{1 -n_{2}}(1 +z(x_{1} ,y_{1}))\vert  \\
 &=\vert d\vert \vert x_{1}^{1 -n_{1}}\vert \vert y_{1}^{1 -n_{2}}\vert \vert 1 +z(x_{1} ,y_{1})\vert  =r^{2 -n}
 \end{align}
 and
 \begin{equation}
	\vert \varepsilon ^{ -i}\widetilde{\beta }_{0} \vert  = \vert \varepsilon ^{ -i}\vert \vert \widetilde{\beta }_{0} \vert  =r^{n/k} .
\end{equation}
If we compose the mapping $(f -W(g^{\frac{1}{k}}) ,g^{\frac{1}{k}}) :T_{r ,r} \rightarrow T_{r^{2 -n} ,r^{n/k}}$ with simple linear mappings of $\mathbb{K}^{2}$ from the left and right: $H_{2} \circ (f -W(g^{\frac{1}{k}}) ,g^{\frac{1}{k}}) \circ H_{1} ,$ where 
\begin{align}H_{1}(x ,y) &=(x_{1}x ,y_{1}y) , \\
H_{2}(x ,y) &=(\frac{x}{dx_{1}^{1 -n_{1}}y_{1}^{1 -n_{2}}} ,\frac{y}{{x_{1}^{{n_1}/{k}}}y_{1}^{n_{2}/k}}) ,\end{align}  we obtain a mapping $T_{1 ,1}$ into $T_{1 ,1}$ of the form
\begin{equation}T_{1 ,1}\! \ni \!(x ,y) \mapsto (x^{1 -n_{1}}y^{1 -n_{2}}(1 +z(x_{1}x ,y_{1}y)) ,x^{\frac{n_{1}}{k}}y^{\frac{n_{2}}{k}}(1 +w(x_{1}x ,y_{1}y)))\! \in \! T_{1 ,1} .
\end{equation}
It suffices to show that the cardinality of the fibres of this mapping is equal to $\vert n_{1} -n_{2}\vert /k .$ Hence we have reduced the problem to finding the number of solutions in $T_{1 ,1}$ of the system 
\begin{equation}
  \left\{\!\!\!\begin{array}{l}
    X^{1 - n_1} Y^{1 - n_2}  (1 + \bar{z} (X, Y)) = \bar{\alpha}\\
    X^{n_1 / k} Y^{n_2 / k}  (1 + \bar{v} (X, Y)) = \bar{\beta}
  \end{array}\right.\!\!\!, \label{37}
\end{equation}
where $\overline{z}(X ,Y) :=z(x_{1}X ,y_{1}Y)  \in \mathbb{K}\mathbb{L}_{\infty }^{X ,Y}$, $\overline{v}(X ,Y) :=v(x_{1}X ,y_{1}Y) \in \mathbb{K}\mathbb{L}_{\infty }^{X ,Y}\!$, $\overline{\alpha } ,\overline{\beta } \in T_{1 ,1}$,  and $\vert \overline{z}(x ,y)\vert  ,\vert \overline{v}(x ,y)\vert  <1$ for $(x ,y) \in T_{1 ,1}$. Notice that  $\overline{z}(X ,Y)$, $\overline{v}(X ,Y)$ have coefficients $\overline{z}_{i ,j}$ and $\overline{v}_{i ,j}$ for which $\left \vert \overline{z}_{i ,j}\right \vert  <1 ,$ $\left \vert \overline{v}_{i ,j}\right \vert  <1.$ This follows from the fact that $z(X ,Y) ,\; v(X ,Y)$ have complex coefficients, $\deg z <0 ,$ $\deg v <0$, and $\left \vert x_{1}\right \vert  =\left \vert y_{1}\right \vert  =r >1.$ 
Even more precisely, from the same fact, we get that $\left \vert \overline{z}_{i ,j}\right \vert$, $\left \vert \overline{v}_{i ,j}\right \vert  \rightarrow 0$ when $i+j \rightarrow -\infty$.

We solve first the ``reduced'' system
\begin{equation}
  \left\{\!\!\!\begin{array}{l}
    X^{1 - n_1} Y^{1 - n_2}  = \bar{\alpha}\\
    X^{n_1 / k} Y^{n_2 / k}   = \bar{\beta}
  \end{array}\right. \label{zz}
\end{equation}
in $T_{1 ,1}$, and next, we prolong these solutions to solutions of (\ref{37}). One can easily check that all the solutions of (\ref{zz}) are $(\overline{x}_{i} ,\overline{y}_{i}) \in T_{1 ,1}$ $(i =1 ,\ldots  ,\left \vert n_{1} -n_{2}\right \vert /k)$, where $\overline{y}_{i}$ are roots of the equation 
$Y^{ (n_{1} -n_{2}) /k} -(\bar{\beta }^{n_{1}-1} \bar{\alpha }^{n_{1}/k}) =0$ and $x_{i} =(\bar{\alpha }\bar{\beta }^{k}/y_{i})$.

In order to prolong these solutions to solutions of (\ref{37}), we apply a version of Hensel's lemma given in \cite[Ch.\ III, {\textsection}4.5, Cor.\ 2]{Bou}, appropriately modified. In notation of this corollary, we put $A : =\{x \in \mathbb{K} :\left \vert x\right \vert  \leq 1\}$ and $\mathfrak{m} : =\{x \in \mathbb{K} :\left \vert x\right \vert  <1\} ,$ $n =2 ,$ $r =0.$ The pair $(A ,\mathfrak{m})$ is a commutative local ring and satisfies Hensel's conditions (given at the beginning of point 5 in Bourbaki's book) because:


1.  $A$ is linearly topologized by the system of ideals $\mathcal{B}_{i} : =\{x \in A :\left \vert x\right \vert  <1/i\}$,

2.  $A$ is separable because $\bigcap _{i}\mathcal{B}_{i} =\{0\}$,

3.  $A$ is complete because $\mathbb{K}$ is,

4. $\mathfrak{m}$ is closed in $A$ because its complement  $\{x \in A :\left \vert x\right \vert  =1\}$ is open,

5. elements of $\mathfrak{m}$ are topologically nilpotent because for any $x \in \mathfrak{m}$ we have $\left \vert x\right \vert  <1.$

Moreover, we take 
\begin{align}f_{1}(X ,Y) &:=X^{1 -n_{1}}Y^{1 -n_{2}}(1 +\overline{z}(X ,Y)) -\bar{\alpha } \\
f_{2}(X ,Y) &:=X^{{n_1}/{k}}Y^{{n_2}/k}(1 +\overline{v}(X ,Y)) -\bar{\beta }. \nonumber
\end{align}
These are not ordinary power series, as is assumed by Bourbaki, but Laurent series at infinity in two variables $X$, $Y$, with coefficients in $A$.  But this case can be reduced to
ordinary power series. It suffices to change the variables with negative
exponents into new ones with opposite (positive) exponents. We will explain
this with an example. For instance, the single equation
\[ X^2 + Y + \frac{1}{Y} + \frac{X}{Y^3} + \frac{1}{X^2 Y^2} + \cdots = 0
   \label{rownanie} \]
involving a Laurent series at infinity can be replaced by the system of three
equations involving a power series
\[ \left\{\begin{array}{l}
     X^2 + Y + V + XV^3 + U^2 V^2 + \cdots = 0\\
     VY -1= 0\\
     UX -1= 0
   \end{array}\right.\!\!\!\!, \label{uklad} \]
where $U$ and $V$ are new variables. The same trick transforms any system of Laurent series at infinity (call it
$S_1$) into a system of power series (call it $S_2$). Clearly, then,
there is a natural 1-1 correspondence between the solutions in $A$ of systems
$S_1$ and $S_2$. Moreover, the jacobian determinant of the new system at any
point $(u_0, v_0, x_0, y_0)$ satisfying $v_0 y_0 = 1$ and $u_0 x_0 = 1$ is
equal to $\pm\, x_0^{- 1} y_0^{- 1}$ times the jacobian determinant of the
original system at the point $(x_0, y_0)$.  Because of these
observations, in the sequel we will continue to proceed with Laurent
series as we would with the ordinary power series.
Thus, note that $(f_1,f_2)$ are \emph{restricted} formal Laurent series, i.e., for each fixed $k$ there exist only finitely many $\overline{z}_{i ,j},$ $\overline{v}_{i ,j}$ not belonging to $\mathcal{B}_{k}$. We also easily check that the jacobian $J_{(f_{1} ,f_{2})}(X ,Y)$ of $(f_{1} ,f_{2})$ is equal to
\begin{equation}J_{(f_{1} ,f_{2})}(X ,Y) =\frac{\left \vert n_{1} -n_{2}\right \vert }{k}X^{1 -n_{1}}Y^{1 -n_{2}}X^{{n_1}/{k}}Y^{{n_2}/k}+t .l .d .
\end{equation}
If we take any solution $(\overline{x}_{i} ,\overline{y}_{i})$ of (\ref{zz}), we obtain\begin{equation}J_{(f_{1} ,f_{2})}(\overline{x}_{i} ,\overline{y}_{i}) =\bar{\alpha }\bar{\beta }\frac{\left \vert n_{1} -n_{2}\right \vert }{k} +R .
\end{equation}We have $\vert \bar{\alpha }\bar{\beta }\frac{\left \vert n_{1} -n_{2}\right \vert }{k}\vert  =1$ and $\left \vert R\right \vert  <1$, i.e., $R \in \mathfrak{m} .$ Hence $J_{(f_{1} ,f_{2})}(\bar{x}_{i} ,\bar{y}_{i})$ is invertible in $A$ and obviously
\begin{equation}(f_{1}(\overline{x}_{i} ,\overline{y}_{i}) ,f_{2}(\overline{x}_{i} ,\overline{y}_{i})) =(\bar{x}_{i}^{1 -n_{1}}\bar{y}_{i}^{1 -n_{2}}\overline{z}(\overline{x}_{i} ,\overline{y}_{i}) ,\bar{x}_{i}^{\,{n_1}/k}\bar{y}_{i}^{\,{n_2}/k}\overline{v}(\overline{x}_{i} ,\overline{y}_{i})) \in \mathfrak{m} \times \mathfrak{m} .
\end{equation}
Hence from the Bourbaki corollary, there is a unique solution $(\overline{\overline{x}}_{i} ,\overline{\overline{y}}_{i})$ of  (\ref{37}) such that $(\overline{\overline{x}}_{i} -\overline{x}_{i},\overline{\overline{y}}_{i} -\overline{y}_{i}) \in \mathfrak{m} \times \mathfrak{m} .$ This gives that system (\ref{37}) has $\left \vert n_{1} -n_{2}\right \vert /k$ solutions.
\end{proof}

Notice that the condition in Proposition \ref{prop2} is fulfilled for $i =0$ because for our fixed $(x_{0} ,y_{0})$ we have $f(x_{0} ,y_{0}) =\alpha _{0}$ and $g^{\frac{1}{k}}(x_{0} ,y_{0}) =\widetilde{\beta }_{0} .$ Hence system (\ref{19}) has exactly $\vert n_{1} -n_{2}\vert /k$ solutions in $T_{r ,r}$ for $i =0.$ Unfortunately, we cannot assert the same for the remaining systems (\ref{19}). Quite the opposite, some of them have no solutions in $T_{r ,r} .$

\begin{proposition}
There exists $i \in \{1 ,\ldots  ,k -1\}$ for which system (\ref{19}) has no solutions in $T_{r ,r} .$

\begin{proof}
From Condition \ref{nu3} there exists a $k$-th root of unity $\varepsilon _{0}$ for which  $f -W(\varepsilon _{0}g^{\frac{1}{k}}) =eX^{n_{1}^{ \prime }}Y^{n_{2}^{ \prime }} +t .l .d. ,\enspace n_{1}^{ \prime } +n_{2}^{ \prime } >2 -n .$ Hence for any $(x ,y) \in T_{r ,r}$ we have $\vert f -W(\varepsilon _{0}g^{\frac{1}{k}})(x ,y)\vert  =n_{1}^{ \prime } +n_{2}^{ \prime } >2 -n .$ Moreover, we have $\vert \alpha _{0}\text{} -W(\widetilde{\beta }_{0})\vert  =r^{2 -n} .$ So, for $i$ such that $\varepsilon ^{i} =\varepsilon _{0}$ the system
\[ \left\{\!\!\!\begin{array}{l}
     (f - W (\varepsilon^i g^{\frac{1}{k}}))  (X, Y) = \alpha_0 - W
     (\widetilde{\beta}_0)\\
     \varepsilon^i g^{1 / k} (X, Y) = \widetilde{\beta}_0
   \end{array}\right. \]
has no solutions in $T_{r ,r} .$ But this system is equivalent to the system (\ref{19}) (for the same $i$). In consequence, for this $i$ system (\ref{19}) has no solutions in $T_{r ,r} .$
\end{proof}

\end{proposition}

From the last proposition, we get that the author's reasoning is wrong and that a proof of the Jacobian Conjecture couldn't be given along his idea.

%

\end{document}